\documentclass[onefignum,onetabnum]{siamart190516}

\usepackage{lipsum}
\usepackage{amsfonts}
\usepackage{graphicx}
\usepackage{overpic}
\usepackage{cleveref}
\usepackage{epstopdf}
\usepackage{mathtools}

\usepackage{algorithm}
\usepackage{algpseudocode}
\algnewcommand{\algorithmicor}{\textbf{ or }}
\algnewcommand{\OR}{\algorithmicor}

\usepackage{xcolor}

\newcommand{\R}{\mathbb{R}}
\renewcommand{\d}{\,\textup{d}}
\renewcommand{\Re}{\textrm{Re}}
\newcommand{\Ra}{\textrm{Ra}}
\renewcommand{\Pr}{\textrm{Pr}}

\newcommand{\omicron}{o}

\ifpdf
  \DeclareGraphicsExtensions{.eps,.pdf,.png,.jpg}
\else
  \DeclareGraphicsExtensions{.eps}
\fi

\newsiamremark{remark}{Remark}
\newsiamremark{example}{Example}
\newsiamremark{hypothesis}{Hypothesis}
\Crefname{hypothesis}{Hypothesis}{Hypotheses}
\newsiamthm{claim}{Claim}

\headers{Optimization of Hopf bifurcation points}{N. Boull\'e, P.~E.~Farrell, and M.~E.~Rognes}

\title{Optimization of Hopf bifurcation points\thanks{Submitted to the editors \today.
\funding{N.B. was supported by the EPSRC Centre for Doctoral Training in Industrially Focused Mathematical Modelling through grant EP/L015803/1 in collaboration with Simula Research Laboratory and an INI-Simons Postdoctoral Research Fellowship. PEF was supported by EPSRC grants EP/R029423/1 and EP/W026163/1. MER has received funding from the European Research Council (ERC) under the European Union's Horizon 2020 research and innovation programme under grant agreement 714892.}}}

\author{Nicolas Boull\'e\thanks{Isaac Newton Institute for Mathematical Sciences, University of Cambridge, CB3 0EH, UK. (\email{nb690@cam.ac.uk})}
\and Patrick E.~Farrell\thanks{Mathematical Institute, University of Oxford, Oxford, OX2 6GG, UK. (\email{patrick.farrell@maths.ox.ac.uk}).}
\and Marie E.~Rognes\thanks{Department of Scientific Computing and Numerical Analysis, Simula Research Laboratory, Oslo, Norway, and Department of Mathematics, University of Bergen, Bergen, Norway. (\email{meg@simula.no}).}}

\usepackage{amsopn}

\begin{document}

\maketitle

\begin{abstract}
We introduce a numerical technique for controlling the location and stability properties of Hopf bifurcations in dynamical systems. The algorithm consists of solving an optimization problem constrained by an extended system of nonlinear partial differential equations that characterizes Hopf bifurcation points. The flexibility and robustness of the method allows us to advance or delay a Hopf bifurcation to a target value of the bifurcation parameter, as well as controlling the oscillation frequency with respect to a parameter of the system or the shape of the domain on which solutions are defined. Numerical applications are presented in systems arising from biology and fluid dynamics, such as the FitzHugh--Nagumo model, Ginzburg--Landau equation, Rayleigh--B\'enard convection problem, and Navier--Stokes equations, where the control of the location and oscillation frequency of periodic solutions is of high interest.
\end{abstract}



\section{Introduction}

Dynamical systems are fundamental in a range of scientific fields including biology, chemistry, physics, medicine, and economics~\cite{strogatz2018nonlinear}. Here, we consider dynamical systems of the form
\begin{equation}
  \label{eq_dyn_syst}
  \frac{\partial u}{\partial t} = F(u,\lambda),
\end{equation}
where $u \in C^1([0,\tau];U)$ is a solution, $\tau$ is the time horizon, $U$ is a suitable Hilbert space of functions defined on a bounded domain $\Omega \subset \R^d$, $d \in \mathbb{N}_+$, $\lambda \in \mathbb{R}$ is a bifurcation parameter, and $F\in C^1(U\times \R;U)$ is a Fr\'echet differentiable operator. Typical examples of operators $F$ are nonlinear partial differential equations (PDEs) such as e.g.~the Navier--Stokes equations. Steady state (or equilibrium) solutions of \cref{eq_dyn_syst} satisfy
\begin{equation}
  \label{eq_steady_system}
  0 \equiv \frac{\partial u}{\partial t} = F(u,\lambda).
\end{equation}
In practice, we employ a finite element method to semi-discretize the time-dependent PDE in space, such that \cref{eq_dyn_syst} can be understood as a system of ODEs with $U=\R^n$ and $F$ mapping $\R^n\times \R$ to $\R^n$.

\begin{figure}[ht!]
\centering
\vspace{0.4cm}
\begin{overpic}[width=\textwidth]{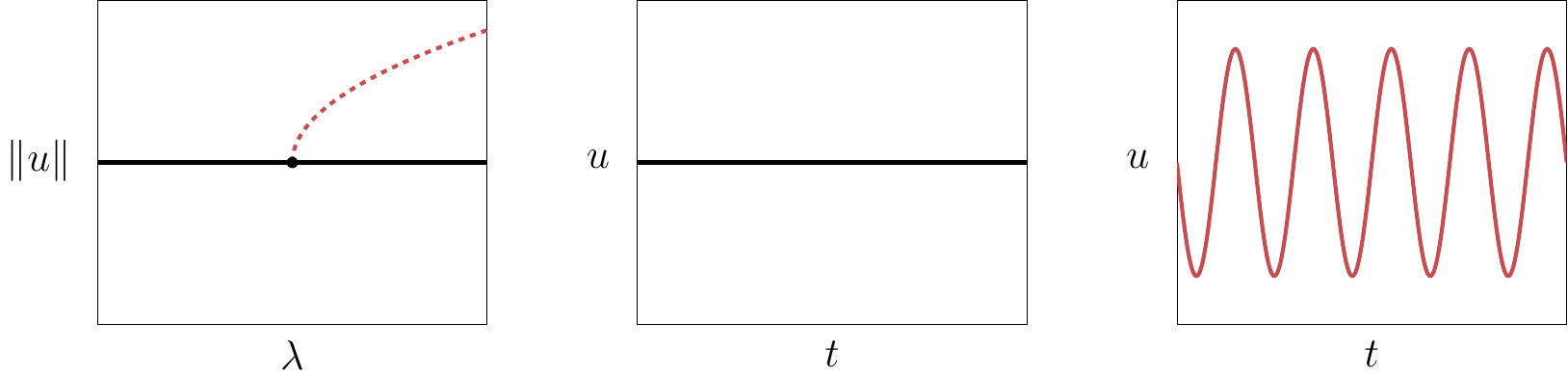}
\put(0,26){(a)}
\put(35,26){(b)}
\put(70,26){(c)}
\end{overpic}
\caption{Illustration of a bifurcation diagram where a branch of transient solutions (dashed line) bifurcates from a steady-state through a Hopf bifurcation (a). The steady-state and periodic solution at the Hopf bifurcation point are respectively depicted in (b) and (c).} 
\label{fig_Hopf}
\end{figure}

The properties of the system \cref{eq_dyn_syst} depend on the value of the bifurcation parameter $\lambda$. For instance, the number of steady-state solutions to \cref{eq_dyn_syst} can vary with the bifurcation parameter $\lambda$ through the birth of branches of solutions $(u, \lambda)$ at specific branching points $(u^\star, \lambda^\star)$ in the bifurcation diagram. Points $(u^\star, \lambda^\star)$ at which steady-state solutions change stability and a periodic solution appears or disappears are known as Hopf bifurcations. By definition, this periodic solution $u$ will satisfy
\begin{equation}
  u(x,t+T) = u(x,t) \quad \text{for } x \in \Omega, t \geq 0
\end{equation}
for some minimal period of oscillation $T>0$. \cref{fig_Hopf} depicts an example of bifurcation diagram of a dynamical system where a branch undergoes a Hopf bifurcation, along with an illustration of the steady-state and periodic solution at the Hopf bifurcation point. Hopf bifurcations are typically of substantial importance in physical and biological systems; in this manuscript we are interested in how and to what extent properties of Hopf bifurcation points can be controlled automatically via numerical optimization.

In this paper, we introduce a numerical method for controlling properties of Hopf bifurcations, such as their location in the bifurcation diagram or the period of the associated periodic solution, by optimizing auxiliary parameters of the model, including the domain $\Omega$. Our algorithm allows us to advance or delay the birth of these periodic solutions using a characterization of Hopf bifurcations as the solutions to an extended system of nonlinear partial differential equations introduced by Griewank and Reddien~\cite{griewank1983calculation}. This extends previous work~\cite{boulle2021control} where we introduced a numerical technique for controlling simple bifurcation structures of physical systems with shape optimization utilizing the Moore--Spence system~\cite{moore1980calculation,seydel1979numerical}.
The proposed algorithm can be applied to many systems of equations modelling sustained oscillations in physical and biological phenomena where Hopf bifurcations naturally arise, such as cardiac cell models~\cite{erhardt2018bifurcation,erhardt2021complex}, the Hodgkin--Huxley model describing action potentials in neurons~\cite{guckenheimer1993bifurcation,hodgkin1952quantitative,wang2005two,xie2008controlling}, the Selkov model of glycolysis~\cite{sel1968self}, the Van der Pol oscillator~\cite{van1927frequency}, spiral waves in reaction-diffusion systems~\cite{barkley1990spiral}, models of cardiac rhythms~\cite{leloup1999limit,roenneberg2008modelling}, the migration of cancer cells~\cite{bruckner2019stochastic}, and the von K\'arm\'an vortex street~\cite{barkley2006linear,duvsek1994numerical,jackson1987finite,provansal1987benard,williamson1996vortex} in fluid dynamics.

The Griewank--Reddien system has been widely used in the literature to accurately localize Hopf bifurcations in various applications such as aeroelasticity, fluid dynamics, and crystals~\cite{chen1997bifurcation,chen1995bifurcation,cliffe2000numerical,dellnitz1990computational,fortin1997localization,luzyanina1996numerical,morton1999hopf,winters1988oscillatory}. However, the majority of the works that aim to control the location and properties of Hopf bifurcations have considered applications to ordinary differential equations (ODEs), with a simple control parameter, and employed methods based on linear stability analysis, center manifold, and/or normal form theory~\cite{murdock2006normal} to find analytical relations~\cite{abed1994stabilization,berns1998feedback,harb2003controlling,jiang2020hopf,li2004hopf,
verduzco2006hopf,wang2020hopf,xie2008controlling,yu2004hopf}. One of the main limitations of these methods is that they are difficult to use to study nonlinear PDEs, especially when the control is more complex than a single real parameter, such as a material parameter field or the shape of the domain. In contrast, our approach can be automatically applied to many examples using a variety of optimization functionals and control parameters.

We have implemented our algorithm in the Firedrake finite element software~\cite{rathgeber2016firedrake}, which interfaces with several libraries such as UFL~\cite{alnaes2014unified}, PETSc~\cite{balay2020petsc}, SLEPc~\cite{hernandez2005slepc}, dolfin-adjoint~\cite{farrell2013automated,mitusch2019dolfin}, ROL~\cite{ridzal2017rapid}, and fireshape~\cite{paganini2021fireshape} for discretizing variational formulations, solving sparse linear systems, solving eigenvalue problems, deriving adjoint models of nonlinear PDEs, solving optimization problems, and performing shape optimization, respectively. The code used to compute the numerical examples in this paper is publicly available on Zenodo for reproducibility purposes~\cite{nicolas_boulle_2021_5235244}.

The paper is organized as follows. We begin in \Cref{sec_charact_hopf} by describing the characterization of Hopf bifurcation points as solutions to an extended system of partial differential equations. We employ the formulation proposed by Griewank and Reddien~\cite{griewank1983calculation}, which generalizes the Moore--Spence system used to characterize simple bifurcations and turning points~\cite{moore1980calculation,seydel1979numerical}. Next, in \Cref{sec_opt}, we describe the algorithm for performing numerical optimization and control of Hopf points with respect to a parameter of the dynamical system or the shape of the domain on which solutions are defined. We present a wide range of numerical examples and applications of the method in \Cref{sec_app} before concluding and discussing further possibilities in \Cref{sec_conc}.

\section{Characterization of Hopf bifurcations}
\label{sec_charact_hopf}

\subsection{Definitions of simple and Hopf bifurcation points}

For a dynamical system of the form~\cref{eq_dyn_syst}, a \emph{simple bifurcation point} is a solution pair $(u^\star, \lambda^\star)$ satisfying \eqref{eq_steady_system} with the property that the number of steady-state solutions to \cref{eq_steady_system} in any neighbourhood of $(u^\star, \lambda^\star)$ changes as $\lambda$ passes $\lambda^\star$. At such points, the Fr\'echet derivative $F_u$ of the operator $F$ with respect to $u$ is non-invertible, with a zero eigenvalue associated with a nonzero eigenfunction $\phi$:
\begin{equation}
  F_u (u^\star, \lambda^\star) [\phi] = 0 .
\end{equation}
Using this property, simple bifurcation points can be characterized using an extended system of equations called the Moore--Spence system~\cite{moore1980calculation,seydel1979numerical}, in which the steady-state $u^\star$, the bifurcation parameter $\lambda^\star$, and the eigenvector $\phi$ are solved for simultaneously.

Following Roose and Hlava\v{c}ek~\cite{roose1985direct}, we assume that there exists a smooth branch of steady-state solutions $(u, \lambda)$, each satisfying \cref{eq_steady_system}, for $\lambda\in(\lambda^\star-\epsilon,\lambda^\star+\epsilon)$ where $\epsilon>0$, and that $u^\star$ is an isolated solution to $F(u, \lambda^\star)=0$ belonging to the branch. Moreover, we assume that $F_u(u^\star, \lambda^\star)$ has a single pair of complex conjugate imaginary eigenvalues $\pm i \mu$ for $\mu \in \mathbb{R}_+$, and with associated eigenfunction $\phi$:
\begin{equation}
  \label{eq_eig_hopf}
  F_u(u^\star, \lambda^\star)[\phi] = i \mu \, \phi .
\end{equation}
Then, under additional regularity conditions on the operator $F$ (cf.~\cite{roose1985direct} and the theoretical results in~\cite{hassard1981theory,marsden2012hopf}), a \emph{Hopf bifurcation} occurs at $\lambda=\lambda^\star$, i.e., a branch of time periodic solutions bifurcates from the branch of steady-state solutions. 

\subsection{The Griewank--Reddien equations}
\label{sec:griewank-reddien}

In this section, we consider the formulation of an extended system of equations to characterize Hopf bifurcations as proposed by Griewank and Reddien~\cite{griewank1983calculation}, which we will next embed in an optimal control setting in \Cref{sec_opt}. Henceforth we consider the semi-discretization in space of \cref{eq_dyn_syst}, as the analysis of Griewank and Reddien considered only ODE problems~\cite{griewank1983calculation,marsden2012hopf}.

While alternative methods~\cite{hassard1981theory,heinemann1981multiplicity,jepson1981numerical} permit the detection of Hopf bifurcations by solving eigenvalue problems during a bifurcation parameter continuation process, the Griewank--Reddien formalism offers a direct approach. One of the key difficulties of the characterization of Hopf points arises from the non-uniqueness of the eigenfunctions, and we will discuss this point before turning to the extended system.

Let $(u^\star,\lambda^\star)$ be a steady-state solution of \cref{eq_steady_system} at a Hopf bifurcation point, and let $\phi_0 = v_0 + i w_0$ be an eigenfunction satisfying \cref{eq_eig_hopf}, with real and imaginary parts $v_0$ and $w_0$ respectively, corresponding to the imaginary eigenvalue $i \mu$ for $\mu > 0$. We observe that $\phi_0$ is only unique up to multiplication by a complex number $z=re^{i\theta}$, and the eigenvalue problem~\cref{eq_eig_hopf} thus admits eigenfunctions of the form 
\begin{equation}
  \phi = re^{i\theta}\phi_0 = r(v_0\cos\theta-w_0\sin\theta) + ir(v_0\sin\theta+w_0\cos\theta),
\end{equation}
for every $r>0$ and $\theta \in [0,2\pi)$. To uniquely determine a solution to \cref{eq_eig_hopf}, we augment the eigenvalue problem by a normalization condition. More precisely, we assume that there exists a suitable \emph{normalization function} $c \in U$ such that at least one of $\langle c,v_0\rangle \neq 0$ or $\langle c,w_0\rangle\neq 0$ holds, where $\langle \cdot, \cdot\rangle$ denotes the inner product of $U$. We aim to fix a solution $\phi = v + i w$, $v, w \in U$ to \cref{eq_eig_hopf} such that
\begin{equation} \label{eq_norm_condition}
\langle c, v \rangle = 0 \qquad \text{and} \qquad \langle c, w \rangle = 1.
\end{equation}
Clearly \cref{eq_norm_condition} is satisfied for $v = r (v_0 \cos \theta - w_0 \sin \theta)$ and $w = r (v_0 \sin \theta + w_0 \sin \theta)$ if $r$ and $\theta$ solve:
\begin{subequations} \label{eq_norm}
\begin{align}
\langle c, v_0\rangle \cos\theta - \langle c, w_0\rangle \sin\theta &= 0,\\
r (\langle c, v_0\rangle \sin\theta + \langle c, w_0\rangle \cos\theta) &= 1.
\end{align}
\end{subequations}
Hence, given $c$ and a $\phi_0$, \cref{eq_norm} provides a normalization procedure to determine $r, \theta$ to ensure that the corresponding eigenfunction $\phi$ satisfies the conditions given by \cref{eq_norm_condition}. 

The Griewank--Reddien formalism combines the eigenvalue problem \cref{eq_eig_hopf}, the real and imaginary components of the steady-state equation \cref{eq_dyn_syst}, as well as the normalization conditions \cref{eq_norm_condition} to read as follows~\cite[Eq.~3.1]{griewank1983calculation}: find the solution field $u \in U$, bifurcation parameter $\lambda \in \R$, frequency $\mu > 0$, and eigenfunction components $v, w \in U$ such that
\begin{equation}
  \label{eq_ext_MS}
  G(u, \lambda, \mu, v,w) =
  \begin{pmatrix}
    F(u, \lambda) \\
    F_u(u,\lambda) [v] + \mu w \\
    F_u(u,\lambda) [w] - \mu v \\
    \langle c, v \rangle  \\
    \langle c, w \rangle -1  
  \end{pmatrix}
  = 0 .
\end{equation}
In the remainder of this paper, we will characterize Hopf bifurcations using the system \cref{eq_ext_MS} and will refer to it as the Griewank--Reddien system. We will also exploit~\eqref{eq_norm} separately to construct better initial guesses for the typically highly nonlinear system.

\begin{remark}
In place of the Griewank--Reddien system~\cref{eq_ext_MS}, other direct formulations could equivalently be used, such as the system proposed by Roose and Hlava{\v{c}}ek~\cite{roose1985direct}:
\begin{equation*}
R(u, \lambda, \mu, \phi) =
\begin{pmatrix}
F(u, \lambda) \\
F_u(u, \lambda)^2[\phi] + \mu^2 \phi\\
\langle c, \phi \rangle \\
\langle \phi, \phi \rangle -1
\end{pmatrix}
= 0 .
\end{equation*}
However, the implementation of this Roose--Hlava{\v{c}}ek system is more challenging in the finite element software Firedrake~\cite{rathgeber2016firedrake}. Alternatively, one could also employ a standard shooting method~\cite{waugh2013matrix} to find periodic solutions by solving a time-dependent system of equations to obtain the solutions and corresponding periods. While this method should allow the control of periodic solutions far from Hopf bifurcations, it requires a transient simulation for each functional evaluation and is therefore much more computationally expensive when employed in an optimization problem.
\end{remark}

\section{Numerical optimization of Hopf bifurcation points}
\label{sec_opt}

A Hopf bifurcation point $(u, \lambda)$ of the dynamical system
\cref{eq_dyn_syst} is associated with various properties, such as the
location $\lambda$ of the bifurcation parameter, the frequency of
oscillation of the emerging periodic solution branch $\mu$, as well as
other properties of the steady-state $u$ at the bifurcation. We target
controlling these properties via an optimal control approach.

\subsection{An optimal control setting for Hopf bifurcations}

Assume that a Hopf bifurcation point $(u, \lambda)$ to the dynamical
system \cref{eq_dyn_syst} with frequency $\mu > 0$ can be expressed as a
function of a control variable $\omicron$ (e.g.~material parameter,
bifurcation location, frequency, domain shape etc.). For a given
objective functional
\begin{equation}
  \mathcal{J} = \mathcal{J} (u, \lambda, \mu) = \mathcal{J} (u, \lambda, \mu) (\omicron),
\end{equation}
we consider the optimization problem constrained by the Griewank--Reddien system~\cref{eq_ext_MS}:
\begin{equation}
  \label{eq_min_problem}
  \min_{\omicron} \, \mathcal{J}(u, \lambda, \mu) \quad \text{subject to} \quad G(u, \lambda, \mu, v, w, \omicron) = 0.
\end{equation}
This general formulation allows for a range of types of optimal control applications such as illustrated by the following four examples.
\begin{example}[Controlling the location of a Hopf bifurcation]
To advance or delay the bifurcation parameter of a Hopf bifurcation point to a target value $\lambda^\star$, we consider the objective functional:
\begin{equation}
  \mathcal{J} = (\lambda - \lambda^\star)^2 / \lambda^{\star 2} .
\end{equation}
\end{example}
\begin{example}[Controlling the frequency of the periodic solution]
  To increase or decrease the frequency of oscillation of the periodic branch arising from the Hopf bifurcation to a target frequency $\mu^\star$, we consider the functional:
  \begin{equation}
    \mathcal{J} = (\mu - \mu^\star)^2/\mu^{\star 2} .
  \end{equation}
\end{example}
\begin{example}[Optimizing a parameter]
To control the Hopf bifurcation with respect to a scalar parameter $a \in \R$ of the PDE, we consider the control variable $o=a$.
\end{example}
\begin{example}[Optimizing the shape of the domain]
To control the Hopf bifurcation with respect to the shape of the domain $\Omega$, we consider the control variable $o=\Omega\in \mathcal{U}_{\textup{ad}}$ where $\mathcal{U}_{\textup{ad}}$ is the set of images of an initial domain under a suitable set of diffeomorphisms~\cite{paganini2021fireshape}.
\end{example}

\subsection{Optimization algorithm}

To minimize the functional $\mathcal{J}$, we introduce an iterative optimization algorithm,
summarized in \cref{alg_branch} and further described below.

\renewcommand{\algorithmicrequire}{\textbf{Input:}}
\renewcommand{\algorithmicensure}{\textbf{Output:}}
\begin{algorithm}
  \caption{Optimization of Hopf bifurcations}
  \label{alg_branch}
  \begin{algorithmic}[1]
    \Require Initial control variable $\tilde{\omicron}$, initial guess $(\tilde{u}, \tilde{\lambda})$ for the Hopf bifurcation point, optimization functional $\mathcal{J}$, normalization function $c$
    \Ensure Optimized control variable $\omicron$
    \State Solve the eigenvalue problem \cref{eq_eig_hopf} around $(\tilde{u}, \tilde{\lambda})$ to generate initial ($\tilde{\mu}, \tilde{v}, \tilde{w})$
    \State Normalize the eigenfunctions $\tilde{v}, \tilde{w}$ by solving \cref{eq_norm} with the given $c$
    \State Solve the system \cref{eq_ext_MS} to obtain an initial solution $(u^{(0)}, \lambda^{(0)}, \mu^{(0)}, v^{(0)}, w^{(0)})$
    \State Initialize optimization step, $k \gets 1$
    \While {termination criteria not satisfied}
    \State Evaluate the objective functional and compute updated control variable $\omicron^{(k)}$
    \State Solve the Griewank--Reddien system to obtain $(u^{(k)}, \lambda^{(k)}, \mu^{(k)}, v^{(k)}, w^{(k)})$
    \If{the regularity conditions are satisfied}
    \State Accept optimization step, $k\gets k+1$
    \Else
    \State Reject optimization step and decrease the step size
    \EndIf
    \EndWhile
  \end{algorithmic}
\end{algorithm}

First, we solve the Griewank--Reddien system~\cref{eq_ext_MS} to locate the Hopf bifurcation to be modified. This system of equations is highly nonlinear and, as a dynamical system may have several Hopf bifurcations, can have multiple solutions. As an example, the Ginzburg--Landau equation example in Section~\ref{sec_GL} has an infinite number of Hopf bifurcation points. We use Newton's method to solve the Griewank--Reddien system, and therefore need a good initial guess to ensure convergence to the target Hopf bifurcation. To this end, we first use deflated continuation~\cite{farrell2016computation,farrell2015deflation} to compute multiple steady-state solutions $u$ to \cref{eq_dyn_syst} by continuation in the bifurcation parameter $\lambda$. We analyze the (linear) stability of each steady-state $u$ found at each parameter $\lambda$ by computing eigenvalues to the Fr\'echet derivative of $F$ at $u$, $F_u(u,\lambda)$, and tracking when one eigenvalue becomes purely imaginary as we increase $\lambda$. We select the steady-state $(\tilde{u}, \tilde{\lambda})$ and eigenfunction $\phi$ with growth rate (real part of the eigenvalue) closest to zero. We define the guess frequency $\tilde{\mu} > 0$ to be the imaginary part of the corresponding eigenvalue, and denote by $v_0$ and $w_0$ the real and imaginary parts of $\phi$, respectively. Last, we define the normalized functions to be
\begin{equation*}
  \tilde{v} = r(v_0 \cos \theta - w_0 \sin \theta), \qquad
  \tilde{w} = r(v_0 \sin \theta + w_0 \sin \theta),
\end{equation*}
where $r>0$ and $\theta\in [0,2\pi)$ are solutions to \cref{eq_norm} in alignment with the discussion in Section~\ref{sec:griewank-reddien} (step 2 in \cref{alg_branch}).

Once a suitable initial guess $(\tilde{u}, \tilde{\lambda}, \tilde{\mu}, \tilde{v}, \tilde{w})$ has been computed, we solve~\cref{eq_ext_MS} to obtain an initial Hopf bifurcation point $(u^{(0)}, \lambda^{(0)}, \mu^{(0)}, v^{(0)}, w^{(0)})$. We employ a trust-region algorithm~\cite{conn2000trust} to solve the optimization problem formulated in \cref{eq_min_problem} and minimize the functional $\mathcal{J}(u, \lambda, \mu)$ (step 6--7 in \cref{alg_branch}). Finally, we check whether certain regularity conditions on the mesh of the domain or the steady-state solution $u^{(k+1)}$ are satisfied to accept or reject the optimization step (cf.~\cite[Sec.~4]{boulle2021control}). In particular, if $u^{(k)}$ is the steady-state solution at the previous step, we reject control updates that do not satisfy the following inequality, to ensure that the optimization remains on the same branch of solutions and does not jump to a secondary Hopf bifurcation:
\[\|u^{(k+1)}-u^{(k)}\|_U\leq C\|u^{(k+1)}\|_U,\]
where $C>0$ is a specified constant determined heuristically. In the examples described in \cref{sec_RB,sec_NS}, we will choose $C=0.2$ and $C=0.05$ to balance the speed of convergence of the optimization with the constraint that we control the correct bifurcation point. If the above inequality is not satisfied, then we reject the optimization step and decrease the trust-region radius. The optimization algorithm terminates when the functional value $\mathcal{J}$ falls below a tolerance $\epsilon>0$, or a maximum number of iterations or a gradient tolerance is reached.

\subsection{Discretization, solvers and software}

The nonlinear (partial) differential equation $F(u,\lambda) = 0$ is discretized and solved using the Firedrake finite element software~\cite{rathgeber2016firedrake} with efficient linear solvers from PETSc~\cite{balay2020petsc}. We solve the eigenvalue problems resulting from the linear stability analysis using the Scalable Library for Eigenvalue Problem Computations (SLEPc)~\cite{hernandez2005slepc}, which is a library interfacing with PETSc for solving large scale eigenvalue problems. In particular, we use the Krylov--Schur algorithm with a shift-and-invert spectral transformation~\cite{kressner2005numerical,stewart2002krylov}.  We leverage the Rapid Optimization Library (ROL)~\cite{ridzal2017rapid} to solve the optimization problem formulated in \cref{eq_min_problem} and minimize the functional $\mathcal{J}(u, \lambda, \mu)$. Inside the optimization algorithm, the update of the control variable is computed using either the dolfin-adjoint library~\cite{farrell2013automated,mitusch2019dolfin}, if the control is a parameter or a function, or the Fireshape optimization toolbox~\cite{paganini2021fireshape} which relies on a moving mesh method~\cite{allaire2006structural,paganini2018higher}, if the control is the shape of the domain. 

\section{Applications}
\label{sec_app}

In this section, we employ the abstract method described in \Cref{sec_opt} for controlling Hopf bifurcations to a wide range of concrete examples. The applications considered include the FitzHugh--Nagumo model simulating the evolution of action potentials in an excitable biological cell such as a neuron or myocyte, the complex Ginzburg--Landau equation used to understand phase transition, as well as the Rayleigh--B\'enard convection problem, and the control of a von K\'arm\'an vortex street described by the Navier--Stokes equations.

\subsection{FitzHugh--Nagumo model}

We first consider a FitzHugh--Nagumo~\cite{fitzhugh1961impulses} model. This dynamical system simplifies the Hodgkin--Huxley model~\cite{hodgkin1952quantitative}, which describes the propagation of action potentials in neurons. This system of ordinary differential equations models the evolution of the transmembrane potential $v$, and a second dimensionless variable $w$ in a cardiac cell, and reads as:
\begin{subequations}
  \label{eq_FN}
\begin{align}
\frac{\partial v}{\partial t} &=c_1v(v-a)(1-v)-c_2 w, \\
\frac{\partial w}{\partial t} &= b(v-c_3 w).
\end{align}
\end{subequations}
Here, $a=-0.12,\,b=0.011,\,c_1=0.15,\,c_2=0.05,\,c_3=0.55$ are given parameters which may be adjusted to model different type of cells~\cite[Sec.~2.4.1]{sundnes2007computing}.
This formulation yields a normalized action potential with a zero resting potential and a peak around $0.9$. Note that a reparametrization of the model can be done to match physiological data~\cite[Sec.~2.4.1]{sundnes2007computing}. 

In this example, we select the parameter $c_1$ as the bifurcation parameter of the model, fix the remaining parameters except $c_2$, and remark that the system transitions to a time-dependent periodic solution to \cref{eq_FN} through a Hopf bifurcation, located at a critical value of the parameter $c_1$. Then, we solve the Griewank--Reddien system~\cref{eq_ext_MS} and obtain a Hopf bifurcation at the critical bifurcation parameter $c_1=0.05$, associated with the steady-state $(v,w)=(0,0)$. Moreover, using linear stability analysis, we find that this state is associated with a pair of imaginary eigenvalues of $\pm i\mu$, where $\mu\approx 2.23\times 10^{-2}$, corresponding to a period of oscillations of $T=2\pi/\mu\approx 277$ ms, and eigenvector $\phi=v_h+iw_h$.

We aim to find the value of the cell parameter $c_2$ such that the action potential duration reflects the action potential duration of a cardiac cell of approximately $T^\star=400$ ms~\cite{boron2016medical, qu2012mechanisms}; i.e., $\mu^\star=1.57\times 10^{-2}$. We formulate this as the following optimization problem:
\begin{equation} \label{eq_MS_FN}
\begin{aligned}
\min_{c_2\in\R} \quad &  \mathcal{J}\coloneqq (\mu-\mu^\star)^2/\mu^{\star 2}\\
\textrm{subject to} \quad & G((v,w),c_2,\mu,v_h,w_h) = 0,
\end{aligned}
\end{equation}
which we solve with the algorithm presented in \cref{sec_opt}, implemented in the Firedrake finite element software~\cite{rathgeber2016firedrake} using the dolfin-adjoint library~\cite{farrell2013automated}. Using a trust-region algorithm implemented in ROL, we are able to minimize the functional $\mathcal{J}$ in \cref{eq_MS_FN} to machine precision in 7 iterations, and obtain an optimized value of $c_2\approx 0.026$. This yields a Hopf bifurcation at $c_1\approx 0.05$ with an associated steady-state of $(v,w)=(0,0)$. We then perform a linear stability analysis to verify that this solution possesses an pair of purely imaginary eigenvalues $\pm i\mu$, with $\mu \approx 1.57\times 10^{-2}$, i.e.~an oscillation period of $T=400$ ms as desired.

\begin{figure}[ht!]
\centering
\vspace{0.4cm}
\begin{overpic}[width=\textwidth]{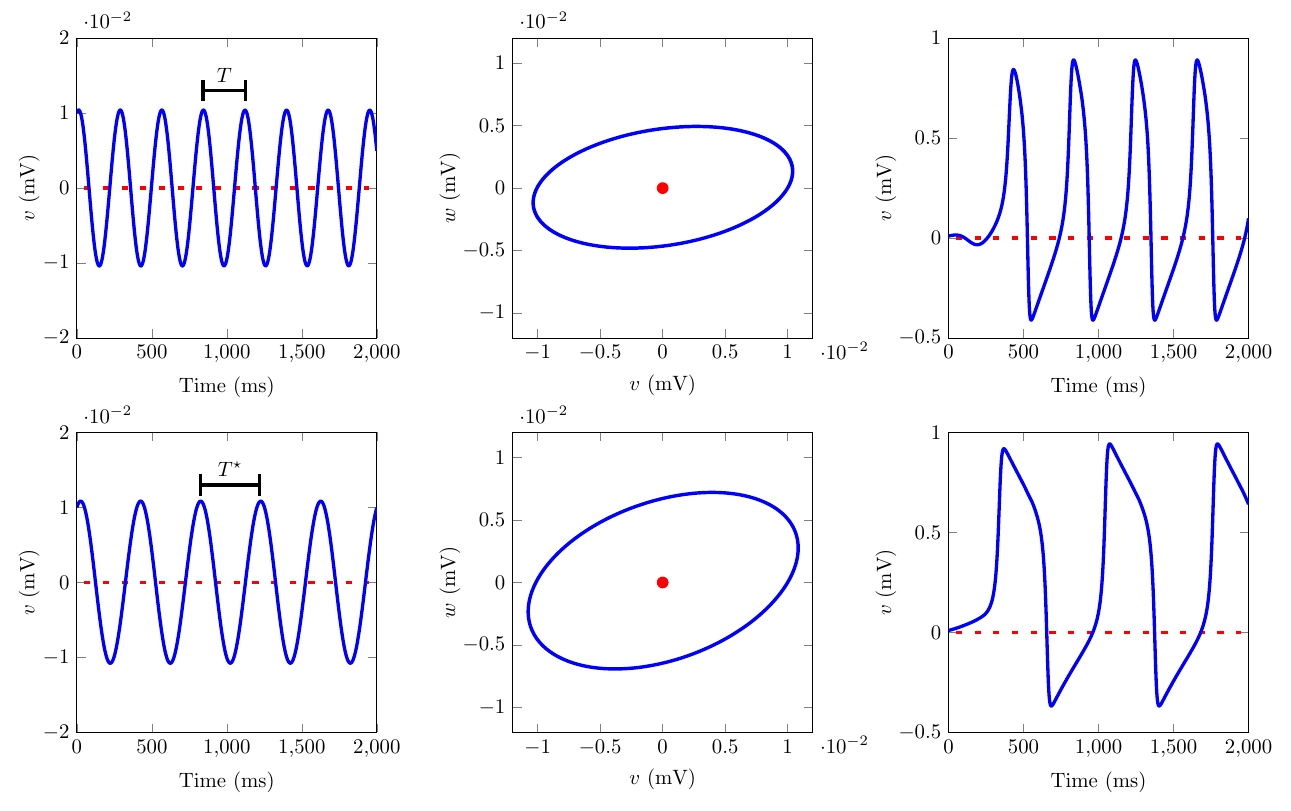}
\put(0,62){(a)}
\put(33,62){(b)}
\put(67,62){(c)}
\put(0,31){(d)}
\put(33,31){(e)}
\put(67,31){(f)}
\end{overpic}
\caption{Initial (a) and optimized (d) periodic potential solutions $v$ to the FitzHugh--Nagumo equations~\cref{eq_MS_FN} at the bifurcation parameter $c_1=0.05$. The dashed red lines show the steady-state solution $(v,w)=(0,0)$ for the respective values of the cell model parameters at the Hopf bifurcation. Panels (b) and (e) display a phase portrait of the variables $v$ and $w$ to illustrate the periodicity of the solution, along with a red dot for the steady-state solution. The periodic solution $v$ at $c_1=0.15$, away from the Hopf bifurcation, is plotted for the initial and optimized parameter $c_2$ in panels (c) and (f), respectively.} 
\label{fig_time_dep_FN}
\end{figure}
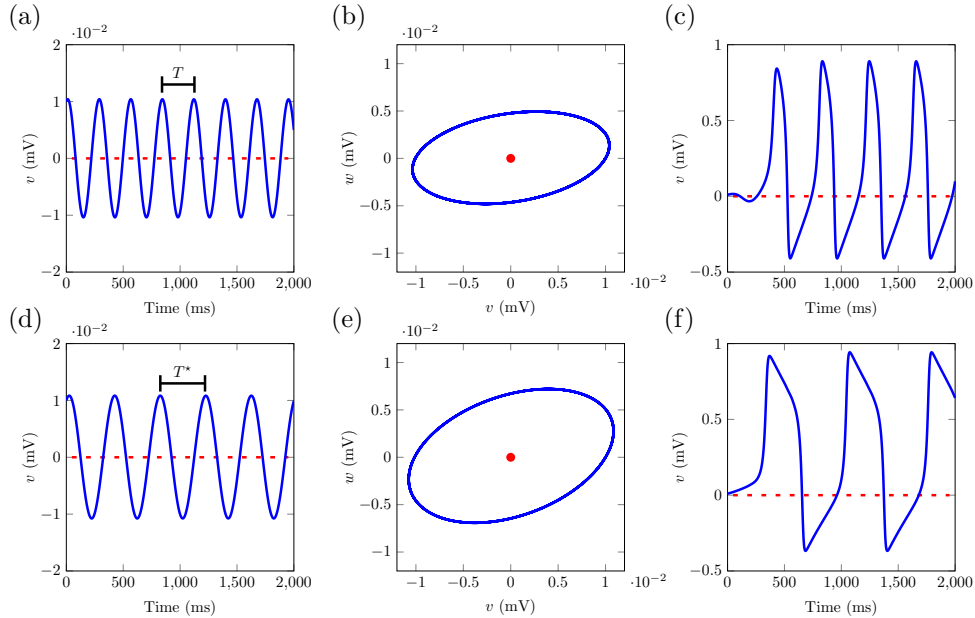

Finally, we employ an explicit Runge-Kutta method of order 5(4)~\cite{dormand1980family} to solve \cref{eq_MS_FN}, starting from the perturbed steady-state $(v,w)=(0.01,0)$, to observe the periodic solution to the FitzHugh--Nagumo equations around the Hopf bifurcation. The resulting periodic potential solutions $v$ for the initial and optimized cell parameters $c_1$ and $c_2$ are displayed in \cref{fig_time_dep_FN}(a) and (d). In the two panels, we observe that the two solutions respectively oscillate with a period of $T=277$ ms and $T=400$ ms, as imposed by the optimization procedure described in this section. Additionally, we report a phase portrait of the variables $v$ and $w$ around the steady-state solution for both cases in \cref{fig_time_dep_FN}(b) and (e). We find in \cref{fig_time_dep_FN}(a) and (d) that the action potential variable $V$ is periodic at the Hopf bifurcation point with sinusoidal oscillations, which are not physiologically realistic. A more relevant solution is obtained by using a higher value of the bifurcation parameter $c_1=0.15$, i.e.,~away from the location of the Hopf bifurcation (see~\cref{fig_time_dep_FN}(c) and (f)). Note that, in this case, the period of the periodic solution has changed and cannot be easily related to the period that we control using the
Griewank--Reddien system. In general, controlling properties of solutions far from the Hopf bifurcation point is considerably more challenging and computationally expensive as it would require the performing of branch continuation within the optimization procedure.

\subsection{Complex Ginzburg-Landau equation} \label{sec_GL}

The Ginzburg--Landau equation is a widely studied nonlinear equation used to describe and understand a wide range of physical phenomena and systems such as phase transitions, nonlinear waves, Bose--Einstein condensates, and liquid crystals~\cite{aranson2002world,bohr1998dynamical,busse1998evolution,crosshohenberg93,kuramoto1984chemical,newell1993order,pismen1999vortices,scott2006encyclopedia}.
This models carries a number of features (dissipation, diffusion, dispersion) that can be controlled using the equation parameters to generate different patterns. From a bifurcation analysis viewpoint, this equation is interesting due to its richness in generating oscillatory or rotating patterns, such as spiral waves~\cite{cross2009pattern,pismen2006patterns}, arising from Hopf bifurcations~\cite{uecker2021numerical}. 

We consider the complex Ginzburg--Landau (CGL) equation with cubic-quintic nonlinearity~\cite{ginzburg1950theory}:
\begin{equation} \label{eq_GL_complex}
\frac{\partial u}{\partial t}=\Delta u+(r+i\nu)u-(c_3+i\mu)|u|^2 u-c_5|u|^4u,\qquad u(x,t)\in\mathbb{C},
\end{equation}
defined on a spatial domain $\Omega\subset\R^2$, with homogeneous Dirichlet boundary conditions. The parameters $\mu$, $\nu$, $c_3$, $c_5$ dictate the dynamics and number of solutions of the equation, and are set by default to $\mu=0.1$, $\nu=1$, $c_3=-1$, $c_5=1$~\cite{uecker2019user}. The parameter $r\geq 0$ plays the role of a bifurcation parameter for the system. 

\begin{figure}[ht!]
\centering
\vspace{0.4cm}
\begin{overpic}[width=\textwidth]{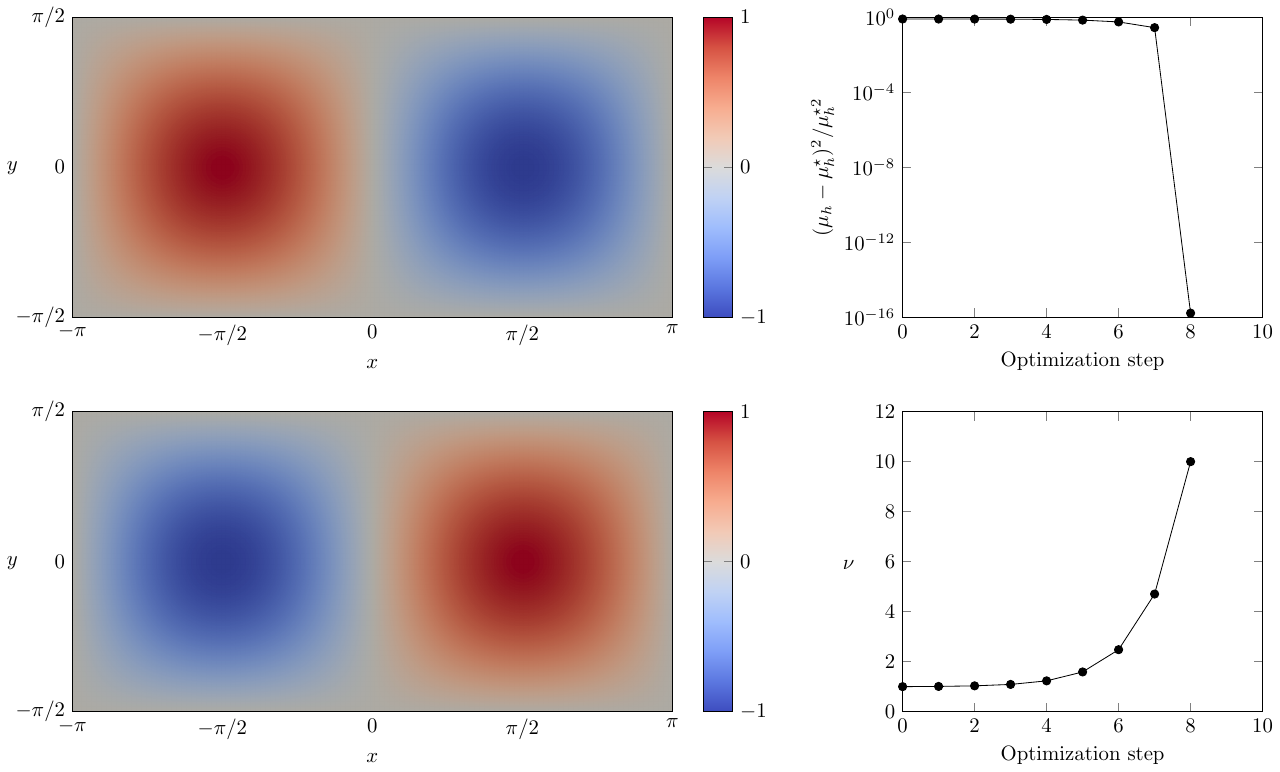}
\put(0,61){(a)}
\put(0,30){(b)}
\put(62,61){(c)}
\put(62,30){(d)}
\end{overpic}
\caption{The normalized two components (a-b) of the real part of the first eigenfunction of the trivial branch $(u_1,u_2)=(0,0)$ at the second Hopf bifurcation located at $r=2$. Panel (c) displays the functional value with respect to the number of optimization steps, while (d) shows the associated value of the control parameter $\nu$.} 
\label{fig_GL_opt}
\end{figure}
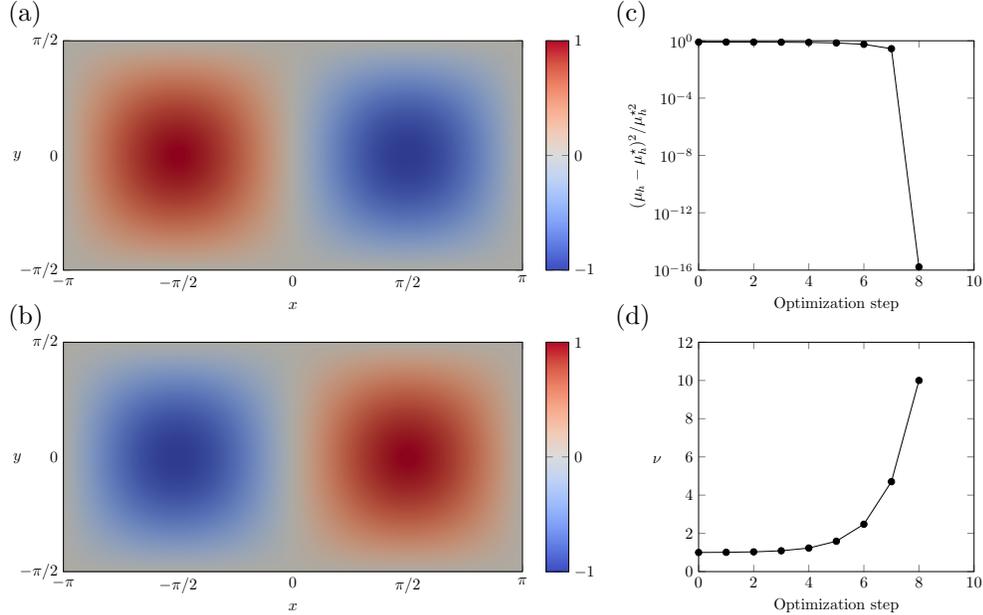

We first decouple \cref{eq_GL_complex} into real and imaginary parts to obtain the following two-component system:
\[\frac{\partial}{\partial t}\begin{pmatrix}
u_1\\ u_2
\end{pmatrix}
=
\begin{pmatrix}
\Delta+r & -\nu \\
\nu & \Delta+r
\end{pmatrix}
\begin{pmatrix}
u_1\\ u_2
\end{pmatrix}
-(u_1^2+u_2^2)
\begin{pmatrix}
c_3 u_1-\mu u_2\\ \mu u_1+c_3 u_2
\end{pmatrix}
-c_5(u_1^2+u_2^2)^2
\begin{pmatrix}
u_1\\ u_2
\end{pmatrix},\]
where $u=u_1+iu_2$. If $\Omega=(-l_1\pi,l_1\pi)\times (-l_2\pi,l_2\pi)$, then it is known that the trivial branch $(u_1,u_2)=(0,0)$ has Hopf bifurcations located at $r=|k|^2\coloneqq k_1^2+k_2^2$, where $k\in \mathbb{Z}/(2 l_1)\times \mathbb{Z}/(2 l_2)$, with an associated pair of imaginary eigenvalues $\pm i\mu_h = \pm i\nu$~\cite{uecker2019user}. The remaining parameters do not influence the location of the Hopf bifurcations with respect to the bifurcation parameter $r$, which depend only on the geometry of the domain $\Omega$ such as its aspect ratio. In this example, we choose the domain $\Omega=(-\pi,\pi)\times (-\pi/2,\pi/2)$, giving the first two Hopf bifurcation points from the trivial branch at $r=5/4$ and $r=2$. To demonstrate the ability of our method to control secondary Hopf bifurcation points in the diagram, provided the initial guess for the Griewank--Reddien system~\cref{eq_ext_MS} is sufficient close to the Hopf point, we focus on the second Hopf bifurcation arising at $r=2$. The two components of the real part of the eigenfunction associated with the pair of imaginary eigenvalues $\pm i$ are illustrated in \cref{fig_GL_opt}(a) and (b). Then, a perturbation of the trivial branch in this direction gives birth to a periodic solution, whose components oscillate between the modes displayed in \cref{fig_GL_opt}(a) and (b) with a period $T = 2\pi$.

We minimize the functional $\mathcal{J}=(\mu_h-\mu_h^\star)^2/\mu_h^{\star 2}$ with respect to the parameter $\nu\in\R$ to control the imaginary part of the eigenvalue associated with the second Hopf bifurcation of the trivial branch as a test case to verify that $\mu_h=\nu$. The variable $\mu_h^\star$ denotes the target value of $\mu_h$ and is set to $\mu_h^\star=10$. The normalization function $c:\Omega\to \R^2$ in the Griewank--Reddien system~\cref{eq_ext_MS} is chosen to be the following function
\begin{equation} \label{eq_norm_const}
c(x,y) = ((x+\pi/2)^2+(y+\pi)^2,-(x+\pi/2)^2-(y+\pi)^2),\qquad x,y\in\Omega.
\end{equation}
With this example, we highlight the importance of the choice of the normalization $c$ to ensure that its inner product with the eigenfunction associated with a purely imaginary eigenvalue is nonzero. In this case, we observe in \cref{fig_GL_opt} that the components of the real part of the eigenfunction satisfy anti-symmetry relations with respect to the $x$ axis. Therefore, the choice of $c=(1,1)$ would be orthogonal to the eigenfunction, which is the reason for the selection of a normalization function breaking the different symmetries of the eigenfunction in \cref{eq_norm_const}. 

We use dolfin-adjoint and ROL to formulate and solve the optimization problem and control the frequency $\mu_h$ at the Hopf bifurcation point with respect to the parameter $\nu$. We report in \cref{fig_GL_opt}(c) the value of the functional throughout the optimization. Our algorithm is able to control the relative frequency to 8 digits of accuracy, corresponding to a minimization of the functional $\mathcal{J}$ to machine precision, in 8 optimization steps. Panel (d) of \cref{fig_GL_opt} displays the corresponding value of the parameter $\nu$ throughout the optimization. As expected, we observe that the target frequency $\mu_h^\star=10$ is reached with the parameter value $\nu=10$, demonstrating the correctness of our method. Note that even in this simple example, where the state $u=(0,0)$ is always a solution to the steady-state equation, the resulting Griewank--Reddien system remains highly nonlinear due to the normalization condition of the eigenfunction and has several solutions. In this case the trivial branch has an infinite number of Hopf bifurcations corresponding to the eigenvalues of the Laplacian on the domain $\Omega$.

\subsection{Rayleigh--B\'enard problem} \label{sec_RB}

We now investigate Hopf bifurcations in a two-dimensional Rayleigh--B\'enard convection problem~\cite{benard1900,benard1927,rayleigh1916}, which models a confined fluid heated from below with a constant temperature difference between the top and bottom of a unit square cell. Several studies have been performed over the past decades to analyze bifurcation structures of the Rayleigh--B\'enard convection problem in various geometries~\cite{bodenschatzetal00,crosshohenberg93,ma2006multiplicity} using numerical methods ranging from arclength continuation and branch-switching techniques~\cite{doedel1981auto,keller1977numerical,uecker2014pde2path}, transient simulations~\cite{boronska2010extreme,boronska2010extreme2}, and deflation~\cite{boulle2021bifurcation,farrell2015deflation}. Additionally, the transition of steady flow structures to oscillatory convection through a Hopf bifurcation has been analyzed numerically and observed in experiments~\cite{colinet1994hopf,ecke1986critical,ecke1992hopf,zhong1991asymmetric} to characterize the associated critical values of the bifurcation parameters. We consider the time-dependent Rayleigh--B\'enard convection problem of an incompressible confined fluid heated from below in a unit square cell domain $\Omega = (0,1)^2$, whose behaviour is governed by the following equations:
\begin{subequations} \label{eq_RB}
\begin{align}
\frac{\partial u}{\partial t} - \Pr\nabla^2u + u\cdot\nabla u  + \nabla p  - \Pr \Ra T\hat{z} &= 0 \quad \text{in }\Omega,\\
\nabla\cdot u&=0 \quad \text{in }\Omega,\\
\frac{\partial T}{\partial t} - \nabla^2 T + u\cdot\nabla T &= 0 \quad \text{in }\Omega,
\end{align}
\end{subequations}
where $u$ is the velocity field, $p$ is the pressure, $T$ is the fluid temperature, $\hat{z}$ is the buoyancy direction, and $\Ra$, $\Pr$ are the Rayleigh and Prandtl numbers. Similarly to~\cite{boulle2021bifurcation}, we assume that the domain has rigid walls with thermally conducting horizontal walls and insulating side walls. That is, we impose the following boundary conditions:
\[u=0 \text{ on }\partial\Omega,\qquad \partial_x T=0\text{ for } x=0,1,\qquad T=1\text{ at }z=0,\qquad T=0\text{ at }z=1.\]
Equations~\cref{eq_RB} are discretized spatially using Taylor--Hood finite elements for the velocity and pressure on triangles (piecewise quadratic and linear polynomials respectively) and piecewise linear polynomials for the temperature using the Firedrake finite element software~\cite{rathgeber2016firedrake}. Note that we solve the steady-state version of~\cref{eq_RB} with $\partial_t u=0$, $\partial_t T=0$, and hence do not need to employ a time-stepping scheme. 

In this example, we fix the Prandtl number to $\Pr=1$ and choose the Rayleigh number as bifurcation parameter. 
We are interested in controlling the location of Hopf bifurcations, i.e.~the critical $\Ra$ at which they arise, in the bifurcation diagram originating from the steady-state solutions to \cref{eq_RB} with respect to the shape of the domain $\Omega$. Following \cref{sec_opt}, we formulate this question as a PDE-constrained shape optimization problem:
\begin{equation} \label{eq_prob_RB}
\begin{aligned}
\min_{\substack{\Omega\in\mathcal{U}_{\textrm{ad}},\, (u,p,T)\in U(\Omega),\\ v,w\in U(\Omega),\, \Ra,\mu>0}} \quad &  \mathcal{J}((u,p,T),\Ra,\mu) \coloneqq (\Ra-\Ra^\star)^2/\Ra^{\star 2}\\
\textrm{subject to} \quad & G((u,p,T),\Ra,\mu,v,w) = 0,
\end{aligned}
\end{equation}
where the system of equations $G$ is defined by \cref{eq_ext_MS} and characterizes the Hopf bifurcation, and $\Ra^\star$ is the target value of the Rayleigh number for the location of the Hopf bifurcation. The set of admissible domains, $\mathcal{U}_{\textrm{ad}}$, in \cref{eq_prob_RB} consists of the image of the initial domain under bi-Lipschitz diffeomorphisms. The PDE-constrained shape optimization algorithm aims to solve \cref{eq_prob_RB} by applying a succession of smooth deformations of the domain~\cite{paganini2021fireshape,paganini2018higher}. The normalization function $c=1$ (in each of the velocity, pressure, and temperature subfunctions) is used to solve \cref{eq_ext_MS}, i.e.~we impose the phase condition $\langle 1,v+iw\rangle=i$ to normalize the eigenfunction $\phi=v+iw$ associated with the eigenvalue $i\mu$. Note that it is important to impose this condition over all velocity, pressure and temperature fields, unlike the Navier--Stokes example presented in \cref{sec_NS}, since the velocity and temperature fields average to zero over the domain $\Omega$ for a symmetric flow.
 
\begin{figure}[ht!]
\centering
\vspace{0.4cm}
\begin{overpic}[width=\textwidth]{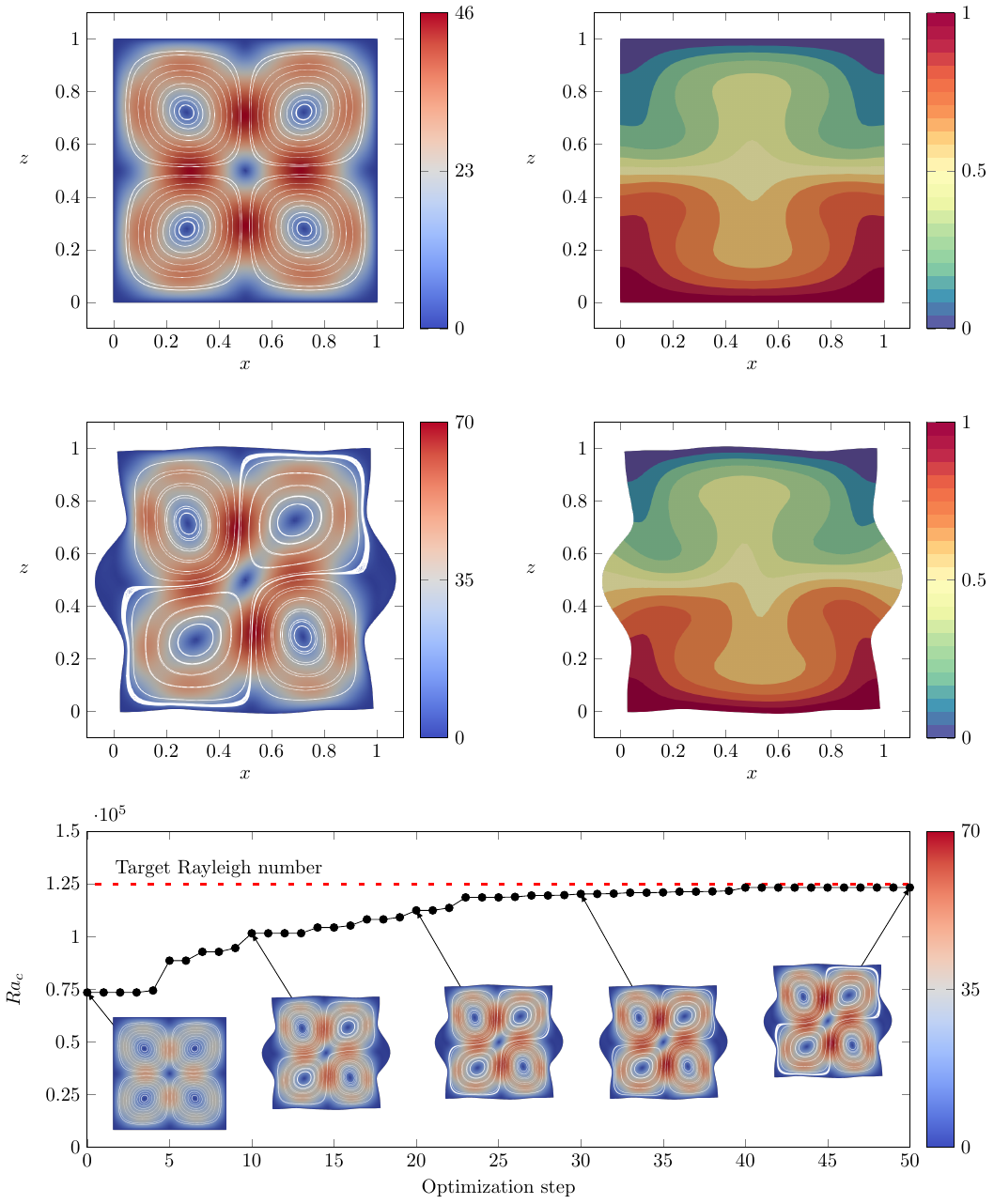}
\put(1,98){(a)}
\put(16.2,101){Velocity}
\put(56,101){Temperature}
\put(1,64){(b)}
\put(1,31){(c)}
\end{overpic}
\caption{Velocity magnitude and streamlines (left) and isotherms (right) of the steady-state solution to \cref{eq_RB} at the Hopf bifurcation located at $\Ra\approx 7.35\times 10^4$ (a). Fluid velocity and temperature flows on the optimized domain at the target Rayleigh number $\Ra^\star=1.25\times 10^5$ (b). Evolution of the critical Rayleigh number, shape of the domain, and corresponding velocity flow at the Hopf bifurcation, during the optimization algorithm (c). Each black dot in the diagram represents the Hopf bifurcation point as the domain is continuously deformed from state (a) to (b).} 
\label{fig_RB_opt}
\end{figure}
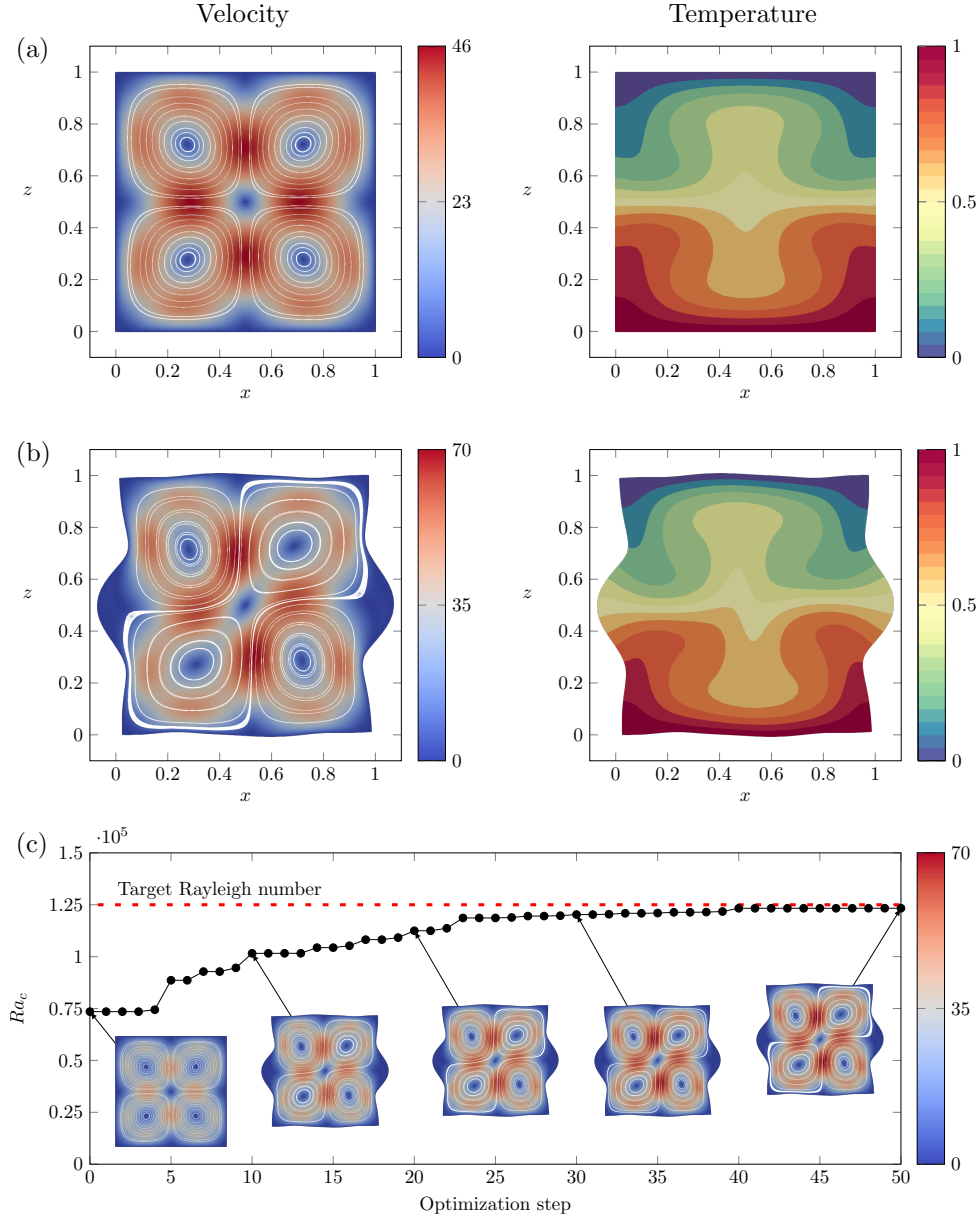

We aim to control the Hopf bifurcation in the branch arising from the 4th bifurcation of the conducting state, located at $\Ra\approx 7.35\times 10^4$, as observed in~\cite[Fig.~8]{boulle2021bifurcation}. The associated flow structure is reminiscent of mechanically coupled stationary convection in bi-layered systems for which a Hopf bifurcation occurs as shown by~\cite[Fig.~1]{colinet1994hopf} and discussed in~\cite{cardin1991nonlinear,cardin1991thermal,nataf1988responsible,rasenat1989theoretical}. We first display the fluid velocity and temperature at the critical Rayleigh number on the initial square domain $\Omega=(0,1)^2$ in \cref{fig_RB_opt}(a). At the Hopf bifurcation, the fluid flow velocity structure is composed of four vortices with alternating rotational directions: clockwise and anticlockwise, together with a symmetric temperature field with respect to the $\mathbb{Z}_2$-symmetry axes of the problem: $x=1/2$ and $z=1/2$. Then, we implement the PDE-constrained optimization problem \cref{eq_prob_RB} in the Fireshape optimization toolbox~\cite{paganini2021fireshape} and solve the resulting optimization problem using a trust-region algorithm~\cite{conn2000trust} implemented in the rapid optimization library (ROL)~\cite{ridzal2017rapid}. We set the target bifurcation parameter to
$\Ra^\star=1.25\times 10^5$ and find the shape deformation of the original domain for which the Hopf bifurcation arises at this value. We report the critical Rayleigh numbers throughout the optimization procedure together with the velocity and temperature fields of the Hopf bifurcation on the optimized domain shape in \cref{fig_RB_opt}(c) and (b), respectively. As we observe in \cref{fig_RB_opt}(c), the optimization algorithm ends by stagnating at a critical Rayleigh number of $\Ra\approx 1.23\times 10^5$, corresponding to an optimization functional value of $\mathcal{J}\approx 2.56\times 10^{-4}$ (cf.~\cref{eq_prob_RB}), i.e.~$1.6\%$ relative error with respect to the target value of the Rayleigh number $\Ra^\star=1.25\times 10^5$. We observe that this error is satisfactory given the coarse mesh discretization of the original domain, and lower errors might be achieved by employing a finer initial mesh. Finally, we observe that the $\mathbb{Z}_2$-symmetries of the temperature field are broken on the optimized solution displayed in \cref{fig_RB_opt}(b) due to the lack of symmetries of the final domain. We do not enforce any constraint to preserve the symmetries of the domain in the optimization formulation~\cref{eq_prob_RB}. Adding further constraints on the shape deformations, such as symmetries or box
constraints, might be of interest for specific applications, as we will see in the next example.

\subsection{Navier--Stokes equations} \label{sec_NS}

\begin{figure}[ht!]
\centering
\vspace{0.2cm}
\begin{overpic}[width=0.8\textwidth]{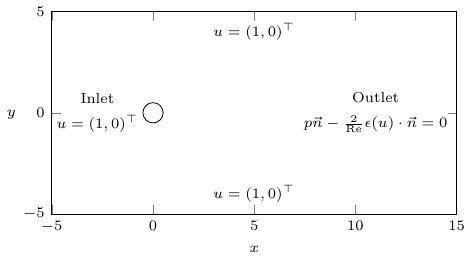}
\end{overpic}
\caption{Initial computational domain $\Omega$ for simulating a von K\'arm\'an vortex street. Here, $\vec{n}$ denotes the unit vector orthogonal to the boundary of the domain $\partial\Omega$. A circular obstacle of radius $0.5$ is located at the coordinates $(x,y)=(0,0)$.} 
\label{fig_NS_domain}
\end{figure}

In this last example, we consider a laminar fluid flow past a circular cylinder in two dimensions. The behaviour of the flow is governed by the Reynolds number $\textrm{Re}$. It is well known that the flow transitions from stationary to periodic at a critical Reynolds number, $\textrm{Re}=\textrm{Re}_c$, through a Hopf bifurcation, and ultimately transitions to turbulence as the Reynolds number increases. The periodic structure is known as a von K\'arm\'an vortex street~\cite{barkley2006linear,duvsek1994numerical,jackson1987finite,provansal1987benard,williamson1996vortex}. The evolution of the fluid flow on the domain $\Omega\subset\R^2$ is modelled by the non-dimensionalized incompressible Navier--Stokes equations:
\begin{subequations} \label{eq_NS}
\begin{align}
\frac{\partial u}{\partial t}-\nabla\cdot \left(\frac{2}{\textrm{Re}}\epsilon(u)\right)+u\cdot\nabla u+\nabla p&=0\quad \text{in }\Omega,\\
\nabla\cdot u&=0\quad \text{in }\Omega,
\end{align}
\end{subequations}
where $\textrm{Re}$ is the Reynolds number, $u$ is the fluid velocity, $\epsilon(u)=\frac{1}{2}(\nabla u+\nabla
u^\top)$, and $p$ is the pressure. We adopt a similar initial domain $\Omega$ as in~\cite{jackson1987finite} and consider a rectangle $(-5,15)\times (-5,5)$ with a circular obstacle centered at the origin of diameter $d=1$. We impose the inflow velocity $u=(1,0)^\top$ at the (left) inlet, top, and bottom of the domain, a no-slip condition on the obstacle boundary, as well as a natural outflow condition at the (right) outlet. We represent the initial computational domain $\Omega$ in \cref{fig_NS_domain}, along with the prescribed boundary conditions. The initial mesh of the domain is generated using Gmsh~\cite{geuzaine2009gmsh} and is composed of $34,444$ triangles, with smaller characteristic length near the obstacle to capture the vortex structures arising in the flow pattern at the critical Reynolds number. Additionally, we impose a symmetric structure of the mesh with respect to the axis $y=0$ to preserve the $\mathbb{Z}_2$-symmetry of the problem. The velocity and pressure are discretized using the Taylor--Hood finite element.

\begin{figure}[ht!]
\centering
\vspace{0.2cm}
\begin{overpic}[width=\textwidth]{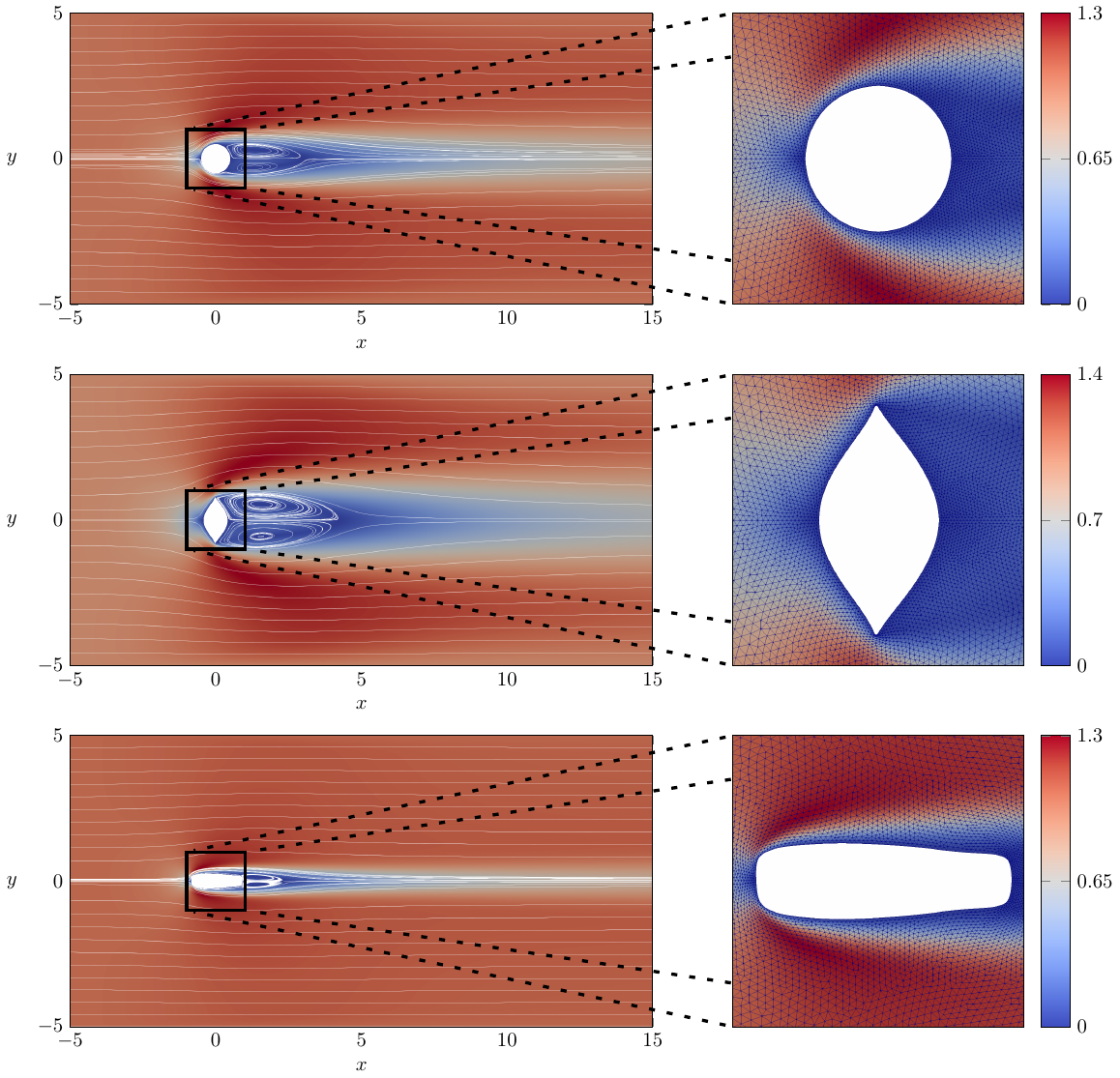}
\put(0,94){(a)}
\put(0,61.5){(b)}
\put(0,29){(c)}
\end{overpic}
\caption{(a) Velocity magnitude and streamlines of the solution at the Hopf bifurcation point located at $\Re^\star\approx 46$, together with a magnification of the mesh around the obstacle at $[-1,1]\times [-1,1]$. (b) Solution to the Navier--Stokes equations on the domain obtained after advancing the Hopf point to $\Re^\star\approx 20$. (c) Similar to (b) with the Hopf bifurcation point delayed to $\Re^\star\approx 200$.} 
\label{fig_NS_karman}
\end{figure}

We first solve the steady-state Navier--Stokes equations at our initial guess for the location of the Hopf bifurcation, $\Re = 46.25$, to obtain a time-independent solution $(u_b,p_b)$. Then, we perform a linear stability analysis using the ansatz $u(x,y,t) = u_b(x,y)+\epsilon\tilde{u}(x,y)e^{\lambda t}$, where $\epsilon\ll 1$, giving the following generalized eigenvalue problem~\cite{barkley2006linear},
\begin{subequations} \label{eq_NS_eig}
\begin{align}
\nabla\cdot \left(\frac{2}{\textrm{Re}}\epsilon(\tilde{u})\right)-u_b\cdot\nabla \tilde{u} - \tilde{u}\cdot\nabla u_b -\nabla \tilde{p}&=\lambda\tilde{u}\quad \text{in }\Omega,\\
\nabla\cdot \tilde{u}&=0\quad \text{in }\Omega,
\end{align}
\end{subequations}
with homogeneous boundary conditions and eigenvalue $\lambda = \sigma+i2\pi f$. Here, $\sigma$ denotes the growth rate and $f$ the frequency of the eigenmode. The steady-state is stable if the largest growth rate is negative and unstable otherwise. We solve \cref{eq_NS_eig} with SLEPc~\cite{hernandez2005slepc} and target the eigenvalues with growth rate $\sigma$ closer to zero. Then, the steady-state solution, leading eigenvalue, and corresponding eigenfunction are used as initial guess for solving the system~\cref{eq_ext_MS} characterizing the Hopf bifurcation. We use the following normalization condition to ensure uniqueness of the solution (see~\cref{eq_norm_condition}),
\[\int_{\Omega}\tilde{u}_x+\tilde{u}_y\,\textup{d} x=i,\]
where $\tilde{u}_x$ and $\tilde{u}_y$ denote the $x$ and $y$ components of the (complex) velocity eigenfunction. We then find a critical Reynolds number of $\textrm{Re}_c\approx 46.23$ with an imaginary pair of eigenvalues of $\lambda=\pm 0.867 i$, i.e., a critical Strouhal number~\cite{strouhal1878besondere,white1999mecanica} of $St_c = f L/U \approx 0.138$, where $L=1$ and $U=1$ are the characteristic length scale and velocity scale. These values are in agreement with the computational study of~\cite{jackson1987finite}, which reported the values $\text{Re}_c=46.136$ and $St_c = 0.13793$ using a similar method for locating the Hopf bifurcation. We report the flow structure (velocity magnitude and streamlines) of the steady-state solution at the Hopf bifurcation in \cref{fig_NS_karman}(a).

We aim to control, i.e.~advance or delay, the Hopf bifurcation by minimizing the functional $\mathcal{J}(\Re) = (\Re-\Re^\star)^2/\Re^{\star 2}$ with respect to the shape of the domain $\Omega$. Here, $\Re^\star$ denotes the target Reynolds number for the location of the Hopf bifurcation and is successively set to $\Re^\star=20$ and $\Re^\star=200$ in the numerical examples presented in this section. In addition, we impose several geometric constraints on the domain. First, we fix the nodes of the mesh at the boundaries of the rectangle $[-5,15]\times [-5,5]$, i.e.~only the inner obstacle may vary. Then, we enforce volumetric and barycentric constraints~\cite{paganini2018higher,schulz2016computational} to ensure that the area and location of the obstacle remain constant throughout the optimization,
\[\int_{\Omega}\d x=\text{constant},\qquad \int_{\Omega}x_1\d x=\text{constant},\qquad \int_{\Omega}x_2\d x=\text{constant},\]
where $x_1$ and $x_2$ denote the two spatial coordinates. These constraints are imposed using an augmented Lagrangian algorithm~\cite[Chapt.~17.3]{nocedal2006numerical} with limited memory BFGS Hessian updates. The subproblems are solved using a trust-region algorithm implemented in ROL.

\begin{figure}[ht!]
\centering
\vspace{0.2cm}
\begin{overpic}[width=\textwidth]{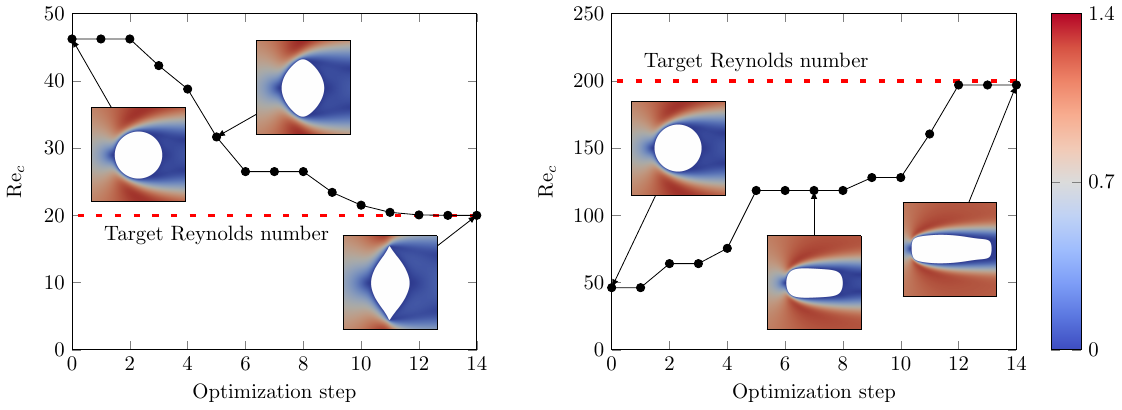}
\end{overpic}
\caption{Evolution of the domain during the optimization algorithm with a target critical Reynolds number of $\text{Re}^\star=20$ (left) and $\text{Re}^\star=200$ (right).  The figures show a magnification of the velocity magnitude profiles on $[-1,1]\times[-1,1]$ around the obstacle, with the same scale and colorbar. The bifurcation parameters achieved at the end of the optimization procedure are respectively equal to $\text{Re}_c = 19.9995$ and $\text{Re}_c = 196.9$.} 
\label{fig_NS_opt}
\end{figure}

We display the original domain and the domains optimized to lead to a Hopf bifurcation at $\Re^\star=20$ and $200$ in \cref{fig_NS_karman}, together with a magnification of the mesh around the obstacles. We first advance the Hopf bifurcation to $\Re^\star=20$ in \cref{fig_NS_karman}(b) and observe that the obstacle is deformed in the vertical direction to reach a final ellipsoid shape with sharp edges. The symmetry of the mesh around the obstacle seems preserved during the optimization. In the left panel of \cref{fig_NS_opt}, we report the evolution of the shape around the obstacle together with the associated critical Reynolds number. The shape optimization procedure successfully converges to a domain with a critical Reynolds number of $\Re_c = 19.9995$. In a second experiment, we aim to delay the birth of instabilities in the fluid flow by finding a shape for which the Hopf bifurcation arises at $\Re^\star=200$. As displayed in the left panel of \cref{fig_NS_opt}, the functional value reaches a plateau at $\Re_c\approx 196.9$, with the shape depicted in \cref{fig_NS_karman}(c). In this case, we observe that an elongated obstacle in the horizontal direction stabilizes the flow for higher Reynolds numbers. This showcases the challenges of the optimization procedure as the deformation of the mesh elements near the obstacle could prevent the functional to decay to machine precision. However, we highlight that we are able to reach a functional value of $\mathcal{J}(\Re) = 2.4\times 10^{-4}$, corresponding to a $1.5\%$ relative error between the target and obtained critical Reynolds number. Several ideas could be implemented to refine these results and reach higher or lower values of the critical Reynolds number, such as preserving the symmetry of the mesh exactly with respect to the axis $y=0$ by defining the problem on the upper half plane, or remeshing during the optimization.

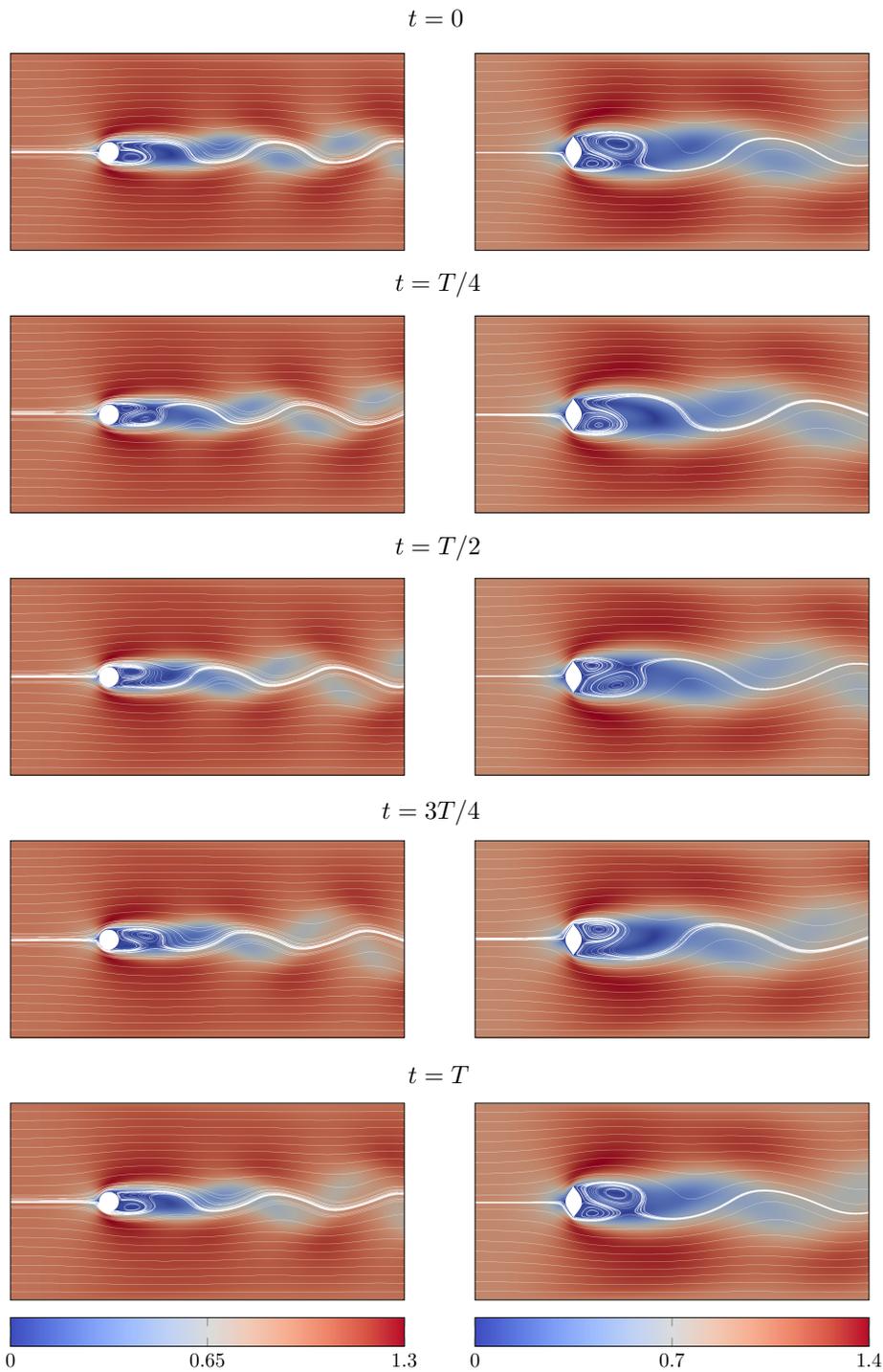
\begin{figure}[htbp]
\centering
\vspace{0.4cm}
\begin{overpic}[width=0.95\textwidth]{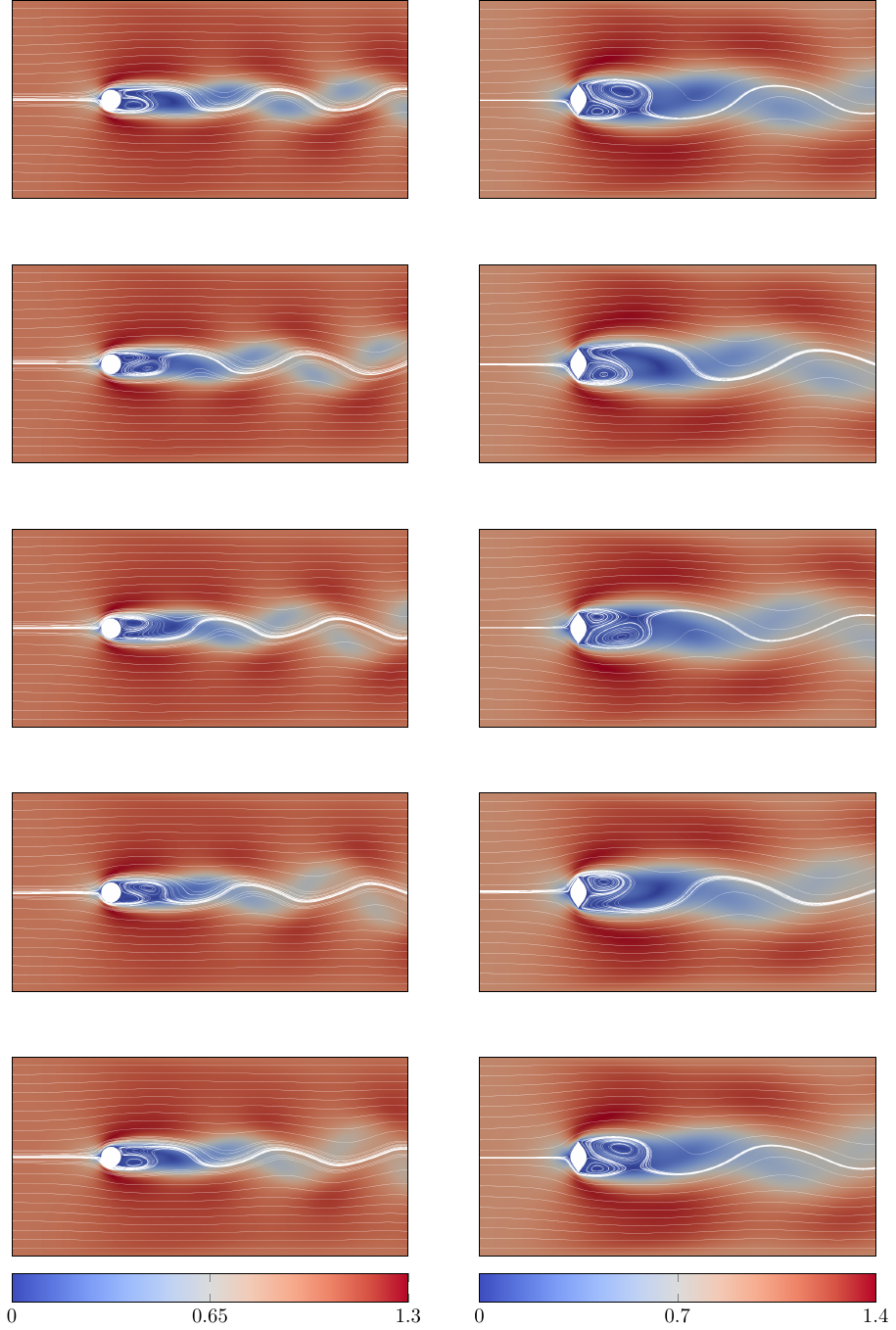}
\put(31,102){$t=0$}
\put(30,82){$t=T/4$}
\put(30,62){$t=T/2$}
\put(29,42){$t=3T/4$}
\put(31,22){$t=T$}
\end{overpic}
\caption{Left: Periodic solution (velocity magnitude and streamlines) to the Navier--Stokes equations, defined on the domain with a cylindrical obstacle, at the Hopf bifurcation $Re\approx 46$ over one time-period of $T\approx 7.14$. Right: Periodic solution on the domain optimized to obtain a Hopf bifurcation at $Re\approx 20$, with a period of $T\approx 9.83$. Movies depicting the evolution of the velocity are available as Supplementary Material.} 
\label{fig_NS_karman_dynamics}
\end{figure}

Finally, we perform time-dependent simulations using both the original domain and the optimized domain depicted in \cref{fig_NS_karman}(b) to observe the von K\'arm\'an vortex street arising at the Hopf bifurcation, at Reynolds numbers $\Re\approx 46$ and $\Re=20$ respectively. We discretize \cref{eq_NS_eig} in time with a Crank--Nicolson time-stepping scheme and use an initial state consisting of a  steady-state to the Navier--Stokes equations perturbed by the eigenmode associated with the Hopf bifurcation, i.e.~$u(x,y,0) = u_b(x,y) +\epsilon\tilde{u}(x,y)$, where $\epsilon$ is chosen such that $\epsilon\|u_b\|_{L^2}/\|\tilde{u}\|_{L^2}=0.05$. We report snapshots of the simulation over one time-period in \cref{fig_NS_karman_dynamics} (movies are available as Supplementary Material) and observe that the velocity profile has a periodic pattern at the expected Reynolds numbers for each simulation.

We have successfully controlled the transition to periodic flow, characterized by a Hopf bifurcation, in two-dimensional laminar flow past a body with respect to the shape of the obstacle. Our approach allows for the systematic manipulation of periodic solutions and may be applied to related problems to analyze the fluid flow past a rotating cylinder~\cite{sierra2020bifurcation} or obstacles with different initial geometries, such as oriented ellipses and triangles~\cite{jackson1987finite}.

\section{Conclusions} \label{sec_conc}

We introduced a robust numerical method for controlling Hopf bifurcations arising in nonlinear dynamical systems. Our algorithm relies on a characterization of Hopf bifurcation points by the Griewank--Reddien system---an extended system of nonlinear partial differential equations, which we embedded into a numerical optimization framework. We applied this procedure to successfully control the location and stability of Hopf bifurcations in several applications, such as the FitzHugh--Nagumo model, the complex Ginzburg--Landau equation, the Rayleigh--B\'enard convection problem, and the Navier--Stokes equations, with respect to a control parameter or the shape of the domain.

\section*{Code availability}
The Firedrake components~\cite{firedrake_zenodo_20211206} and code~\cite{nicolas_boulle_2021_5235244} used to produce the numerical examples presented in this paper are available on Zenodo. The code is also distributed on GitHub at \url{https://github.com/NBoulle/Hopf\_Control}.

\section*{Acknowledgments}
We thank Karoline J{\ae}ger and Aslak Tveito for discussions on the FitzHugh--Nagumo model.

\bibliographystyle{siamplain}
\bibliography{biblio}

\end{document}